\documentclass[a4paper]{article}     
\usepackage{theorem}                 
\usepackage{amssymb}
\usepackage{epsfig}
\newtheorem{theorem}{Theorem}
\newtheorem{corollary}{Corollary}
\newtheorem{proposition}{Proposition}

\newcommand{\CC}{{\mathbb{C}}}
\newcommand{\HH}{{\mathbb{H}}}
\newcommand{\PP}{{\mathbb{P}}}

\newcommand{\RR}{{\mathbb{R}}}
\newcommand{\ZZ}{{\mathbb{Z}}}

\newcommand{\NN}{{\mathbb{N}}}

\newcommand{\eps}{\varepsilon}

\begin{document}

\title{The Poincar\'{e} series of some special quasihomogeneous surface singularities}
\author{Wolfgang Ebeling
\thanks{Partially supported by the DFG-programme ''Global methods in complex
geometry'' (Eb 102/4--1); 2000 Mathematics Subject Classification: Primary 14J17,
32S25, 32S40, 13D40; Secondary 14J28, 11H56}
\\
\parbox{9cm}{\small
 \begin{center} Institut f\"{u}r Mathematik, Universit\"{a}t Hannover, \\ 
        Postfach 6009, D-30060 Hannover, Germany \\
        E-mail: ebeling@math.uni-hannover.de
\end{center}}
}
\date{}

\maketitle

\begin{abstract} 
In \cite{Ebeling01} a relation is proved between the Poincar\'{e} series of the
coordinate algebra of a two-dimensional quasihomogeneous isolated hypersurface 
singularity and the
characteristic polynomial of its monodromy operator. We study this relation
for Fuchsian singularities and show that it is connected with the mirror symmetry of
K3 surfaces and with automorphisms of the Leech lattice. We also indicate relations
between other singularities and Conway's group.
\end{abstract}

\section*{Introduction}

K.~Saito \cite{Saito94, Saito98} has introduced a duality
between polynomials which are products of cyclotomic
polynomials. He has shown that V.~I.~Arnold's strange duality
between the 14 exceptional unimodal hypersurface singularities
is related to such a duality between the characteristic
polynomials of the monodromy operators of the singularities.
Moreover, he has observed that the dual polynomials pair
together to the characteristic polynomial of an automorphism of
the Leech lattice. It is now well-known that Arnold's strange
duality is related to the mirror symmetry of $K3$ surfaces (see
e.g.\
\cite{Dolgachev96}).

The author \cite{Ebeling96, Ebeling98} has shown that these features still hold in a
certain way for the extension of Arnold's strange duality discovered by C.~T.~C.~Wall
and the author \cite{EW85}.

The 14 exceptional singularities and the singularities involved in the extension of
Arnold's strange duality are examples of Fuchsian singularities. By this we mean the
following. Let $\Gamma\subset {\rm PSL}(2,\RR)$ be a finitely generated Fuchsian group of
the first kind. Let $A_k$ denote the $\CC$-vector space of 
$\Gamma$-automorphic forms of weight $2k$, $k \geq 0$, and let 
$A = \bigoplus_{k=0}^\infty A_k$
be the algebra of $\Gamma$-automorphic forms. Then $(X,x):=(\mbox{Spec} \, A, \frak{m})$
where $\frak{m}:= \bigoplus_{k=1}^\infty A_k$ is a normal surface singularity which is
called a Fuchsian singularity. 

In this paper we show that the observation of K.~Saito extends in a certain sense to
all Fuchsian singularities which are isolated complete intersection singularities
(abbreviated ICIS in the sequel), although for some of the polynomials there is a
singularity missing. 

The basis of this duality forms a relation between the Poincar\'{e} series of the
coordinate algebra of such a singularity and the
characteristic polynomial of its monodromy operator which was considered in
\cite{Ebeling01}. There we introduced polynomials $\phi_A(t)$ and $\psi_A(t)$ as follows. 
Let
$(X,x)$ be a normal surface singularity  with good
$\CC^\ast$-action. The coordinate algebra $A$ is a graded algebra. We consider the
Poincar\'{e} series $p_A(t)$ of $A$. Let $\{ g; b; (\alpha_1, \beta_1), \ldots , (\alpha_r, \beta_r)
\}$ be the orbit invariants of $(X,x)$. 
We define
\begin{eqnarray*}
\psi_A(t) & := & (1-t)^{2-r}\prod_{i=1}^r (1-t^{\alpha_i}), \\
\phi_A(t) & := & p_A(t)\psi_A(t).
\end{eqnarray*}

If $(X,x)$ is a Fuchsian singularity, then we show that the
polynomial $\phi_A(t)$ can be interpreted
as the characteristic polynomial of a certain Coxeter
element $c_\infty$ (Proposition~\ref{Poincare}). In the case of Arnold's strange duality,
the Coxeter element $c_\infty$ is the monodromy operator of
the dual singularity.

If $(X,x)$ is a Fuchsian ICIS, then
we derive from the results of \cite{Ebeling01} that we can slightly modify the
polynomial
$\phi_A(t)$ to a rational function $\tilde{\phi}_A(t)$
such that the dual (in Saito's sense) of $\tilde{\phi}_A(t)$ is the
characteristic polynomial $\phi_M(t)$ of the monodromy operator of $(X,x)$ (or a slightly
modified polynomial $\phi^\flat_M(t)$). 
We relate this duality to the mirror symmetry
of K3 surfaces. 

We show that for a Fuchsian ICIS the rational
function $\tilde{\phi}_A(t)$ and its dual pair together to the characteristic
polynomial of an automorphism of the Leech lattice. Moreover, 
we consider  
the quasihomogeneous hypersurface singularities in
$\CC^3$ with Milnor number
$\mu=24$. It was already observed by K.~Saito that the polynomials $\phi_M(t)$
of these singularities are self-dual and are also characteristic
polynomials of automorphisms of the
Leech lattice. We show that the same is true for the polynomials $\phi^\flat_M(t)$ of
some ICIS in $\CC^4$ with $\mu=25$. Finally we indicate 5 constructions which lead
from singularities to self-dual characteristic polynomials of
automorphisms of the Leech lattice and we show that all such polynomials can be
obtained in a suitable way.

The paper is organized as follows. In Sect.~1 we recall the structure of normal
surface singularities with good $\CC^\ast$-action and the definition of Fuchsian
singularities. In Sect.~2 we review the relevant results of \cite{Ebeling01} about the
relation between Poincar\'{e} series and monodromy. In Sect.~3
we consider the polynomial $\phi_A(t)$ of a Fuchsian
singularity and show that it is the characteristic polynomial of a Coxeter element
$c_\infty$. In Sect.~4 we derive the duality among the Fuchsian ICIS and relate it to
the mirror symmetry of K3 surfaces. Finally we discuss the relation to automorphisms of
the Leech lattice in Sect.~5.

The author thanks R.-O.~Buchweitz for pointing out an error in an earlier version of the
paper and C.~T.~C.~Wall for pointing out that the list of ICIS in
$\CC^4$ with $\mu=25$ in that version was incomplete.

\section{Quasihomogeneous surface singularities} \label{QSS}

Let $(X,x)$ be a normal surface singularity with a good $\CC^\ast$-action. So $X$ is a
normal two-dimensional affine algebraic variety over
$\CC$ which is smooth outside its {\em vertex} $x$. Its coordinate ring $A$ has the
structure of a graded $\CC$-algebra 
$A = \bigoplus_{k=0}^\infty A_k$, $A_0=\CC$, and $x$ is defined by the maximal ideal
$\frak{m}= \bigoplus_{k=1}^\infty A_k$. 

According to I.~Dolgachev \cite{Dolgachev75}, there exist a simply connected Riemann
surface $\cal D$, 
a discrete cocompact subgroup $\Gamma$ of $\mbox{Aut}({\cal D})$ and a line bundle
$\cal L$ on $\cal D$ to which the action of $\Gamma$ lifts such that 
$$A_k = H^0({\cal D}, {\cal L}^k)^\Gamma.$$

Let $Z:={\cal D}/\Gamma$. By \cite[Theorem~5.1]{Pinkham77a} (see also
\cite[Theorem~5.4.1]{Wagreich81}), there exist a divisor
$D_0$ on $Z$, $p_1, \ldots , p_r \in Z$, and integers $\alpha_i$, $\beta_i$
with $0<\beta_i < \alpha_i$ and $(\alpha_i,\beta_i)=1$ for 
$i=1, \ldots, r$ such that 
$$A_k =L \left( kD_0 + \sum_{i=1}^r \left[ k \frac{\alpha_i - \beta_i}{\alpha_i} \right]
p_i \right).$$
Here $[x]$ denotes the largest integer $\leq x$, and $L(D)$ for a divisor $D$ on $Z$
denotes the linear space of meromorphic functions $f$ on $Z$ such that $(f) \geq -D$.
We number the points $p_i$ so that $\alpha_1
\leq \alpha_2 \leq \ldots \leq \alpha_r$. Let $g$ be the genus of $Z$ and define $b :=
{\rm degree}\, D_0 +r$. Then $\{ g; b; (\alpha_1, \beta_1), \ldots , (\alpha_r, \beta_r)
\}$ are called the {\em orbit invariants} of $(X,x)$, cf.\ e.g.\ \cite{Wagreich83}.
Define
$\mbox{vdeg}({\cal L}):= -b + \sum_{i=1}^r \frac{\beta_i}{\alpha_i}$.

Now assume that $(X,x)$ is Gorenstein. By \cite{Dolgachev83}, there exists an integer
$R$ such that ${\cal L}^{-R}$ and the tangent bundle $T_{\cal D}$ of $\cal D$ are
isomorphic as $\Gamma$-bundles and
\begin{eqnarray*}
R \cdot \mbox{vdeg}({\cal L}) & = & 2 -2g - r + \sum_{i=1}^r \frac{1}{\alpha_i}, \\
R\beta_i & \equiv & 1 \ \mbox{mod} \, \alpha_i, \quad i=1, \ldots , r.
\end{eqnarray*}
Following \cite[3.3.15]{Dolgachev??} we call $R$ the {\em exponent} of $(X,x)$. 
Since $b$ and the $\beta_i$ are
determined by the $\alpha_i$ and the number $R$, we write the
orbit invariants also as $g; \alpha_1, \ldots, \alpha_r$.

If $R=1$, then we have the following situation. In this case $\Gamma \subset PSL(2,\RR)$ is
a finitely generated cocompact {\em Fuchsian group of the first kind}. This means that
$\Gamma$ acts properly discontinuously on $\HH$ and that the quotient $Z =\HH / \Gamma$
is a compact Riemann surface. The divisor $D_0$ is the canonical divisor, the points 
$p_1, \ldots, p_r \in Z$ are the branch points of the map $\HH \to Z$, 
 $\alpha_i$ is the
ramification index over
$p_i$, and $\beta_i=1$ for $i=1, \ldots ,r$. Hence the orbit invariants are
$$\{ g; 2g-2+r; (\alpha_1,1), \ldots, (\alpha_r,1) \}.$$ 
We follow \cite{Looijenga84} in calling $(X,x)$  a {\em Fuchsian singularity}. The orbit
invariants $\{g ; \alpha_1, \ldots , \alpha_r\}$ are also called the {\em signature}
of $\Gamma$. 

\section{The Poincar\'{e} series}
Let $(X,x)$ be a normal surface singularity  with good
$\CC^\ast$-action with orbit invariants $\{ g; b; (\alpha_1, \beta_1), \ldots , (\alpha_r, \beta_r)
\}$.
Let $p_A(t)$ be the Poincar\'{e} series of the coordinate algebra $A$ of $(X,x)$.
We define
\begin{eqnarray*}
\psi_A(t) & := & (1-t)^{2-r}\prod_{i=1}^r (1-t^{\alpha_i}), \\
\phi_A(t) & := & p_A(t)\psi_A(t), \\
\end{eqnarray*}

%

%

Let $(X,x)$ be an ICIS
with weights $q_1, \ldots, q_n$ and degrees $d_1, \ldots , d_{n-2}$. Then its
Poincar\'{e} series is given by (see e.g.\ \cite[Proposition~(2.2.2)]{Wagreich83})
$$p_A(t) = \frac{\prod_{i=1}^{n-2} (1-t^{d_i})}{\prod_{j=1}^n (1-t^{q_j})}.$$
Therefore $p_A(t)$, $\psi_A(t)$
and
$\phi_A(t)$ are rational functions of the form
$$\phi(t)= \prod_{m|h} (1 - t^m)^{\chi_m} \quad \mbox{for} \
\chi_m \in \ZZ \mbox{ and for some } h \in \NN.$$

Given a rational function
$$\phi(t)=\prod_{m|h} (1 - t^m)^{\chi_m},$$
K.~Saito \cite{Saito94} has defined a dual rational function 
$$\phi^\ast(t) = \prod_{k | h} (1-t^k)^{-\chi_{h/k}}.$$

In \cite{Ebeling01} we proved the following results. 
For integers $a_1, \ldots , a_r$ we denote by $\langle a_1, \ldots ,a_r
\rangle$  their least common multiple and by
$(a_1, \ldots , a_r)$ their greatest common divisor.

\begin{theorem} \label{phi_M3}
Let $(X,x)$ be a quasihomogeneous hypersurface singularity in $\CC^3$. Consider the
rational function
$$\tilde{\phi}_A(t):= \frac{\phi_A(t)}{(1-t)^{2g}}.$$
Then $\tilde{\phi}_A^\ast(t)$ is the characteristic polynomial of the classical
monodromy operator of $(X,x)$.
\end{theorem}

\begin{theorem} \label{phi_M4a}
Let $(X,x)$ be a quasihomogeneous ICIS in $\CC^4$ with weights $q_1$, $q_2$, $q_3$, $q_4$
and  degrees $d_1$, $d_2$.  Assume that $g(z_1,z_2,z_3,z_4)=z_1z_4+z_2z_3$. Define  
$$\tilde{\phi}_A(t):=
\frac{\phi_A(t)(1-t^{d_2})}{(1-t)^{2g}(1-t^{d_1})}, \quad \phi_M^\flat(t):=
\frac{\phi_M(t)}{(1-t)}.$$ 
Then we have  $\tilde{\phi}_A^\ast(t)=\phi_M^\flat(t)$.
\end{theorem}

\begin{theorem} \label{phi_M4b}
Let $(X,x)$ be a quasihomogeneous ICIS in $\CC^4$ with weights $q_1$, $q_2$, $q_3$, $q_4$
and  degrees $d_1$, $d_2$.  Assume that either

(A) $g(z_1,z_2,z_3,z_4)=z_1^q+z_2z_3$
and $f(z_1,z_2,z_3,z_4)=f'(z_1,z_2,z_3)+z_4^p$ for some integers $p,q \geq 2$ where
$q|d_2$, or

(B) $g(z_1,z_2,z_3,z_4)=z_1^q+(z_2-z_3)z_4$
and $f(z_1,z_2,z_3,z_4)=az_1^q+z_2(z_3-z_4)$ for some $a \in \CC$, $a \neq 0,1$, and some
integer $q \geq 2$ and $p:=2$. 

Define  
\begin{eqnarray*}
\tilde{\phi}_A(t) & := &
\frac{\phi_A(t)(1-t^{d_2})^{p-1}(1-t^{\frac{d_1}{q}})(1-t^{\frac{d_2}{p}})}
{(1-t)^{2g}(1-t^{d_1})(1-t^{\frac{d_2}{q}})^p}, \\ 
\phi_M^\flat(t) & := & \frac{\phi_M(t)(1-t^q)^p}{(1-t)^{p-1}(1-t^{\langle p,q
\rangle})^{(p,q)}}. 
\end{eqnarray*}
Then we have
$\tilde{\phi}_A^\ast(t)=\phi_M^\flat(t)$.
\end{theorem}

\section{The Poincar\'{e} series of a Fuchsian singularity}
Now let $(X,x)$ be a Fuchsian singularity.
Let $p_A(t)$ be the Poincar\'{e} series of the algebra $A$. We have
$$p_A(t)= \frac{1+(g-2)t+(g-2)t^2+t^3}{(1-t)^2} + \sum_{i=1}^r
\frac{t^2(1-t^{\alpha_i-1})}{(1-t)^2(1-t^{\alpha_i})}.$$

We shall now show that the polynomial $\phi_A(t)$ can also be
interpreted as the characteristic polynomial of a certain operator. 

Let $(X,x)$ be a normal surface singularity with good $\CC^\ast$-action. Then $X$ can
be compactified to $\bar{X}$ in a natural way (see \cite{Pinkham77a}). 
The variety $\bar{X}$ has $r$ cyclic quotient singularities of type
$(\alpha_1,\alpha_1-\beta_1), \ldots , (\alpha_r,\alpha_r-\beta_r)$
along $\bar{X}_\infty:=\bar{X}-X$ \cite[Lemma~4.1]{Pinkham77a}.

\begin{figure}\centering
\unitlength1cm
\begin{picture}(8.5,6.5)
\put(0.5,0){\includegraphics{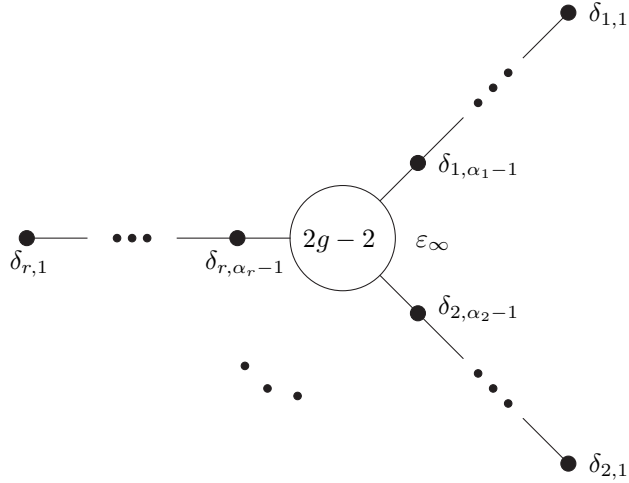}}
\put(8.1,6){$\delta_{1,1}$}
\put(8.1,0){$\delta_{2,1}$}
\put(6.1,4){$\delta_{1,\alpha_1-1}$}
\put(6.1,2.1){$\delta_{2,\alpha_2-1}$}
\put(4.3,3){$2g-2$}
\put(5.8,3){$\eps_\infty$}
\put(3,2.7){$\delta_{r,\alpha_r-1}$}
\put(0.4,2.7){$\delta_{r,1}$}
\end{picture}
\caption{Dual graph of $\tilde{X}_\infty$}
\label{Fig1}
\end{figure}

Now assume that $(X,x)$ is Fuchsian. Then $\beta_i=1$ for all $i=1, \ldots , r$. A
cyclic quotient singularity of type $(\alpha,\alpha-1)$ is a singular point of type
$A_{\alpha-1}$. Let $\pi: \tilde{X} \to \bar{X}$ be the minimal resolution of the
singularities of $\bar{X}$ along $\bar{X}_\infty$. 
The preimage $\tilde{X}_\infty$ of $\bar{X}_\infty$
consists of the strict transform $\eps_\infty$ of $\bar{X}_\infty$ and $r$ chains
$\delta_{i,1}, \ldots, \delta_{i,\alpha_i-1}$, $i=1, \ldots , r$, of rational curves of
self-intersection $-2$ which intersect according to the dual graph shown in
Fig.~\ref{Fig1}. By the adjunction formula, the self-intersection number of the curve
$\eps_\infty$ is $2g-2$.

Let $M_\infty$ be the abstract lattice spanned by these curves, i.e., the free
$\ZZ$-module spanned by $\delta_{1,1}, \ldots, \delta_{1,\alpha_1-1}; \ldots ;
\delta_{r,1}, \ldots , \delta_{r,\alpha_r-1}; \eps_\infty$ with the bilinear form
$\langle \ , \
\rangle$ given by the intersection numbers. Let $U$ be a unimodular hyperbolic plane,
i.e., a free $\ZZ$-module  with basis
$\{f_1, f_2\}$ satisfying $\langle f_1, f_1 \rangle = \langle f_2, f_2 \rangle = 0$,
$\langle f_1, f_2 \rangle = \langle f_2 , f_1\rangle =1$. We shall consider the lattice
$N = M_\infty \oplus U$. Special automorphisms of this lattice are defined as follows.
If $\delta \in N$ is a vector of squared length $\langle \delta, \delta \rangle = \pm
2$, then 
$$s_\delta (x) = x - \frac{2\langle x, \delta \rangle}{\langle \delta, \delta \rangle}
x \quad \mbox{for } x \in N$$
defines the {\rm reflection} $s_\delta$ corresponding to $\delta$. For $f \in U$ with
$\langle f,f \rangle =0$ and $w
\in M_\infty$  we define a transformation
$\psi_{f,w}$, called the {\em Eichler-Siegel transformation} corresponding to $f$ and
$w$, by the formula
$$\psi_{f,w}(x) = x + \langle x, f \rangle w - \langle x,w \rangle f - \frac{1}{2}
\langle w,w \rangle \langle x, f \rangle f$$
for $x\in N$ (cf.~\cite{Ebeling87}).

\begin{proposition} \label{Poincare}
Let $(X,x)$ be a Fuchsian singularity. Then
the polynomial $\phi_A(t)$ is the characteristic polynomial of the operator
$$c_\infty = s_{\delta_{1,1}} \cdots s_{\delta_{1,\alpha_1-1}} \cdots 
s_{\delta_{r,1}}  \cdots  s_{\delta_{r, \alpha_r-1}}\psi_{f_1,\eps_\infty} s_{f_1-f_2}.$$
\end{proposition}

\noindent {\em Proof.} 
In order to compute the characteristic polynomial of $c_\infty$, we want to apply
\cite[Ch.~V, \S ~6, Exercice~3]{Bourbaki68}. For two vectors $u,v \in N$ the {\em
pseudo-reflection} $s_{u,v}$ is defined by
$$s_{u,v}(x) = x - \langle x, v \rangle u.$$
An easy calculation shows that the Eichler-Siegel transformation $\psi_{f,
w}$ can be written as a
product of two pseudo-reflections as follows:
$$\psi_{f,w} = s_{\tilde{w}, f} s_{f,w} \quad \mbox{where }
\tilde{w} := \frac{1}{2}\langle w, w \rangle f -w.$$
In the basis
$$
\eps_1  := \frac{1}{2}\langle \eps_\infty, \eps_\infty \rangle f_1 - \eps_\infty,\quad
\eps_2  :=  f_1 , \quad
\eps_3  :=  f_1-f_2,
$$
the operator $\psi_{f_1,\eps_\infty} s_{f_1-f_2}$ can be written as $s_1 s_2 s_3$ where
$$s_i(\eps_j) = \eps_j - a_{ij} \eps_i, \quad 1 \leq i,j \leq 3, \quad
(a_{ij})_{1 \leq j \leq 3}^{1 \leq i \leq 3} = \left( \begin{array}{ccc} 0&0&-1\\
-\langle \eps_\infty, \eps_\infty \rangle & 0 & 0 \\
\frac{1}{2}\langle \eps_\infty, \eps_\infty \rangle & 1 & 2 \end{array} \right).$$
By \cite[loc.cit.]{Bourbaki68} we get using $\langle \eps_\infty ,
\eps_\infty \rangle = 2g-2$
$$
\det (tI-s_1s_2s_3) =  \left| \begin{array}{ccc} t-1 & 0 & -t \\
                                                   2-2g & t-1 & 0 \\
                                                    g-1 & 1 & t+1 
                                 \end{array} \right| 
  = 1 + (g-2)t + (g-2)t^2 + t^3.
$$
This proves Proposition~\ref{Poincare} for the case $r=0$. The general case also
follows by using the formula of \cite[loc.cit.]{Bourbaki68} for the determinant of the
matrix $tI - c_\infty$ and the Laplace expansion formula.

\addvspace{3mm}

\begin{figure}\centering
\unitlength1cm
\begin{picture}(8.5,6.5)
\put(0.5,0){\includegraphics{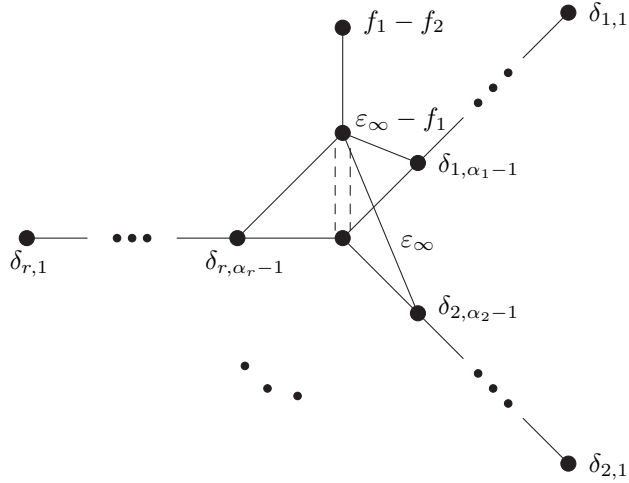}}
\put(8.1,6){$\delta_{1,1}$}
\put(8.1,0){$\delta_{2,1}$}
\put(6.1,4){$\delta_{1,\alpha_1-1}$}
\put(6.1,2.1){$\delta_{2,\alpha_2-1}$}
\put(3,2.7){$\delta_{r,\alpha_r-1}$}
\put(0.4,2.7){$\delta_{r,1}$}
\put(5.6,3){$\eps_\infty$}
\put(5,4.6){$\eps_\infty-f_1$}
\put(5.1,5.9){$f_1-f_2$}
\end{picture}
\caption{The graph $\tilde{\theta}_{\alpha_1, \ldots , \alpha_r}$}
\label{Fig2}
\end{figure}

If $g=0$, then we have $\langle \eps_\infty , \eps_\infty \rangle = -2$ and from
\cite[Sect.~5.1, (c1)]{Ebeling87} we conclude that
$$ \psi_{f_1, \eps_\infty} = s_{\eps_\infty} \circ s_{\eps_\infty -f_1}. $$
Therefore the operator $c_\infty$ coincides with the Coxeter element corresponding to
the graph of Fig.~\ref{Fig2} and its characteristic polynomial is already given
in \cite[Sect.~3.4]{Ebeling87} (unfortunately with a misprint). 

\addvspace{3mm}

\noindent {\bf Remark 1} Consider the polynomial 
$$\psi_A(t)= (1-t)^{2-r} \prod_{i=1}^r (1-t^{\alpha_i})$$
for integers $r \geq 1$, $\alpha_i \geq 1$, $i=1, \ldots , r$.
By \cite[p.~98]{Ebeling87}, this polynomial is the characteristic polynomial of the
Coxeter element corresponding to the graph $\theta_{\alpha_1, \ldots , \alpha_r}$ which
is obtained from the graph of Fig.~\ref{Fig2} by omitting the vertex corresponding to
the vector $f_1-f_2$. The dual rational function is
$$\psi_A^\ast(t) = \frac{(1-t^{\langle \alpha_1, \ldots , \alpha_rÊ\rangle})^{r-2}}
{\prod_{i=1}^r (1-t^{\langle \alpha_1, \ldots , \alpha_r \rangle / \alpha_i})}.$$
For $r \geq 3$ this rational function 
is the Poincar\'{e} series of the Brieskorn-Hamm ICIS $(V_B(\alpha_1, \ldots ,
\alpha_r),0)$ where
$$V_B(\alpha_1, \ldots, \alpha_r):=\{ z \in \CC^r \, | \, b_{i1}z_1^{\alpha_1}+
\ldots + b_{ir}z_r^{\alpha_r}=0; \ i=1, \ldots , r-2 \}$$
and $B=(b_{ij})$ is a sufficiently general $(r-2) \times r$-matrix of complex
numbers. If $(V_B(\alpha_1, \ldots , \alpha_r),0)$ is a simply elliptic singularity,
then the graph $\theta_{\alpha_1, \ldots , \alpha_r}$ is a Coxeter-Dynkin diagram of
this singularity
\cite{Ebeling87}. Therefore we again obtain the identity $\phi_M(t)=p^\ast_A(t)$
(respectively $\phi^\flat_M(t)= p^\ast_A(t)$) in this case (cf.\ \cite{Ebeling01}).

\section{Fuchsian ICIS and mirror symmetry of K3 surfaces} \label{MS}
We now consider Fuchsian ICIS. We first indicate the classification of Fuchsian ICIS. 

Let $(X,x)$ be a Fuchsian singularity.
If the algebra
$A$ is generated by 3 elements, then
$(X,x)$ is a hypersurface singularity in $\CC^3$. 
These cases were classified by I.~Dolgachev \cite{Dolgachev74}, I.~G.~Sherbak 
\cite{Sherbak78}, and Ph.~Wagreich
\cite{Wagreich80}. If $A$ is generated by 4 elements, then
one has an ICIS in $\CC^4$. These cases were classified by Ph.~Wagreich 
\cite{Wagreich81} (see also \cite{Wagreich83}). In the case $g=0$,
$(X,x)$ is a minimally elliptic singularity
\cite[Proposition~5.5.1]{Wagreich81} and the classification can also be derived from
H.~Laufer's results \cite{Laufer77}. There are a few more cases of ICIS of higher
embedding dimension. More precisely we have:

\begin{theorem} \label{ICIS}
Let $(X,x)$ be a Fuchsian singularity with signature $\{g ;  \alpha_1, \ldots , 
\alpha_r\}$. There is an ICIS with this signature if and only if the signature
satisfies one of the following conditions:

{\rm (i)} $g=0$, $r=3$, $\alpha_1=2$, $\alpha_2=3$, and $7 \leq \alpha_3 \leq 10$;

{\rm (ii)} $g=0$, $r=3$, $\alpha_1=2$, $\alpha_2 \geq 4$, 
and $9 \leq \alpha_2+\alpha_3 \leq 12$;

{\rm (iii)} $g=0$, $r=3$, $\alpha_1, \alpha_2, \alpha_3 \geq 3$, 
and $10 \leq \alpha_1+\alpha_2+\alpha_3 \leq 13$;

{\rm (iv)} $g=0$, $r \geq 4$ and $9 \leq \alpha_1+\alpha_2+ \ldots +\alpha_r \leq 12$;

{\rm (v)} $g=1$, $r\geq 1$, $\alpha_1+\alpha_2+ \ldots +\alpha_r \leq r+4$;

{\rm (vi)} $2 \leq g \leq 4$, $r \geq 0$, $\alpha_1+\alpha_2+ \ldots +\alpha_r \leq
r+4-g$;

{\rm (vii)} $g=5$, $r=0$.
\end{theorem}

\noindent {\em Proof.} The algebra $A$ is generated by 3 or 4 elements if and only if
we have one of the cases (i)--(vi) and $Z$ is non-hyperelliptic in the
cases $3;2$ and $4;$. In the remaining cases one can easily show that there is an ICIS with these
invariants, using \cite[Theorem~3.3]{Wagreich81}. Conversely, we show
that in the other cases there are no ICIS. The embedding dimension is given in
\cite[loc.cit.]{Wagreich81}. If $g=0$, then the singularity
$(X,x)$ is minimally elliptic and the result follows from \cite[Theorem~3.13]{Laufer77}.
In the remaining cases we can apply \cite[Lemma~3.9]{VanDyke88} to show that $(X,x)$ is
not an ICIS if the signature does not satisfy the conditions of Theorem~\ref{ICIS}.
This proves Theorem~\ref{ICIS}.

\addvspace{3mm}

\begin{table}\centering
\caption{Fuchsian ICIS with $g=0$ and $r \leq 4$} \label{table1}
\begin{tabular}{|c|c|c|c|c|}
\hline
$g; \alpha_1, \ldots , \alpha_r$ & Weights & Name & Equation(s) &  $\mu$ \\
\hline
$0;2,3,7$ & 6,14,21/42 & $E_{12}$ & $x^7+y^3+z^2$ & 12 \\
$0;2,3,8$ & 6,8,15/30 & $Z_{11}$ & $x^5+xy^3+z^2$ & 11 \\
$0;2,3,9$ & 6,8,9/24  & $Q_{10}$ & $x^4+xz^2+y^3$ & 10 \\
$0;2,3,10$ & 6,8,9,10/16,18 & $J'_9$ & $\left\{ \begin{array}{c} xw+y^2 \\ x^3+yw+z^2
\end{array} \right\}$  & 9 \\
$0;2,4,5$ & 4,10,15/30 & $E_{13}$ & $x^5y+y^3+z^2$ & 13 \\
$0;2,4,6$ & 4,6,11/22 & $Z_{12}$ & $x^4y+xy^3+z^2$ & 12 \\
$0;2,4,7$ & 4,6,7/18 & $Q_{11}$ & $x^3y+xz^2+y^3$ & 11 \\
$0;2,4,8$ & 4,6,7,8/12,14 & $J'_{10}$ & $\left\{ \begin{array}{c} xw+y^2 \\ x^2y+yw+z^2
\end{array} \right\}$  & 10 \\
$0;2,5,5$ & 4,5,10/20 & $W_{12}$ & $x^5+y^4+z^2$ & 12 \\
$0;2,5,6$ & 4,5,6/16 & $S_{11}$ & $x^4+xz^2+y^2z$ & 11 \\
$0;2,5,7$ & 4,5,6,7/11,12 & $L_{10}$ & $\left\{ \begin{array}{c} xw+yz \\ x^3+yw+z^2
\end{array} \right\}$  & 10 \\
$0;2,6,6$ & 4,5,6,6/10,12 & $K'_{10}$ & $\left\{ \begin{array}{c} xw+y^2 \\ x^3+z^2+w^2
\end{array} \right\}$  & 10 \\
$0;3,3,4$ & 3,8,12/24 & $E_{14}$ & $x^8+y^3+z^2$ & 14 \\
$0;3,3,5$ & 3,5,9/18 & $Z_{13}$ & $x^6+xy^3+z^2$ & 13 \\
$0;3,3,6$ & 3,5,6/15 & $Q_{12}$ & $x^5+xz^2+y^3$ & 12 \\
$0;3,3,7$ & 3,5,6,7/10,12 & $J'_{11}$ & $\left\{ \begin{array}{c} xw+y^2 \\ x^4+yw+z^2
\end{array} \right\}$  & 11 \\
$0;3,4,4$ & 3,4,8/16 & $W_{13}$ & $x^4y+y^4+z^2$ & 13 \\
$0;3,4,5$ & 3,4,5/13 & $S_{12}$ & $x^3y+xz^2+y^2z$ & 12 \\
$0;3,4,6$ & 3,4,5,6/9,10 & $L_{11}$ & $\left\{ \begin{array}{c} xw+yz \\ x^2y+yw+z^2
\end{array} \right\}$  & 11 \\
$0;3,5,5$ & 3,4,5,5/8,10 & $K'_{11}$ & $\left\{ \begin{array}{c} xw+y^2 \\ x^2y+z^2+w^2
\end{array} \right\}$  & 11 \\
$0;4,4,4$ & 3,4,4/12 & $U_{12}$ & $x^4+y^3+z^3$ & 12 \\
$0;4,4,5$ & 3,4,4,5/8,9 & $M_{11}$ & $\left\{ \begin{array}{c} xw+yz \\ x^3+(y+z)w
\end{array} \right\}$  & 11 \\
$0;2,2,2,3$ & 2,6,9/18 & $J_{3,0}$ & $x^9+y^3+z^2$ & 16 \\
$0;2,2,2,4$ & 2,4,7/14 & $Z_{1,0}$ & $x^7+xy^3+z^2$ & 15 \\
$0;2,2,2,5$ & 2,4,5/12 & $Q_{2,0}$ & $x^6+xz^2+y^3$ & 14 \\
$0;2,2,2,6$ & 2,4,5,6/8,10 & $J'_{2,0}$ & $\left\{ \begin{array}{c} xw+y^2 \\
x^5+xy^2+yw+z^2               \end{array} \right\}$ & 13 \\
$0;2,2,3,3$ & 2,3,6/12 & $W_{1,0}$ & $x^6+y^4+z^2$ & 15 \\
$0;2,2,3,4$ & 2,3,4/10 & $S_{1,0}$ & $x^5+xz^2+y^2z$ & 14 \\
$0;2,2,3,5$ & 2,3,4,5/7,8 & $L_{1,0}$ & $\left\{ \begin{array}{c} xw+yz \\
x^4+xy^2+yw+z^2               \end{array} \right\}$ & 13 \\
$0;2,2,4,4$ & 2,3,4,4/6,8 & $K'_{1,0}$ & $\left\{ \begin{array}{c} xw+y^2 \\
x^4+xy^2+z^2+w^2               \end{array} \right\}$ & 13 \\
$0;2,3,3,3$ & 2,3,3/9 & $U_{1,0}$ & $x^3y+y^3+z^3$ & 14 \\
$0;2,3,3,4$ & 2,3,3,4/6,7 & $M_{1,0}$ & $\left\{ \begin{array}{c} xw+yz \\
(x^2+w)(y+z)               \end{array} \right\}$ & 13 \\
$0;3,3,3,3$ & 2,3,3,3/6,6 & $I_{1,0}$ & $\left\{ \begin{array}{c} x^3+(y-z)w \\
ax^3+y(z-w)               \end{array} \right\}$ & 13 \\
\hline 
\end{tabular}
\end{table}

\begin{table}\centering
\caption{Fuchsian ICIS with $g=0$ and $r \geq 5$ or $g > 0$} \label{table2}
{\tabcolsep4pt
\begin{tabular}{|c|c|c|c|c|}
\hline
$g; \alpha_1, \ldots , \alpha_r$ & Weights  & Name & Equation(s) & $\mu$   \\
\hline
$0;2,2,2,2,2$ & 2,2,5/10  & $NA^1_{0,0}$ & $x^5+y^5+z^2$ & $16$ \\ 
$0;2,2,2,2,3$ & 2,2,3/8  & $VNA^1_{0,0}$ & $x^4+y^4+yz^2$ & $15$ \\ 
$0;2,2,2,2,4$ & 2,2,3,4/6,6 & $\alpha1^{(1)}$ & $\left\{ \begin{array}{c} x^3+yw \\
 xw+y^3+z^2 \end{array} \right\}$ & $14$ \\
$0;2,2,2,3,3$ & 2,2,3,3/5,6 & $\alpha1^{(2)}$ & $\left\{ \begin{array}{c} xw+yz \\
x^3+y^3+zw \end{array} \right\}$ & $14$ \\
$0;2,2,2,2,2,2$ & 2,2,2,3/4,6 & $\delta1$ & $\left\{ \begin{array}{c} xy+z^2 \\
x^3+y^3+z^3+w^2 \end{array} \right\}$ & $15$ \\
$1;2$ & 1,4,6/12 & $J_{4,0}$ & $x^{12}+y^3+z^2$ & $22$  \\ 
$1;3$ & 1,3,5/10 & $Z_{2,0}$ & $x^{10}+xy^3+z^2$ & $21$ \\ 
$1;4$ & 1,3,4/9 & $Q_{3,0}$ & $x^9+xz^2+y^3$ & $20$ \\ 
$1;5$ & 1,3,4,5/6,8 & $J'_{3,0}$ & $\left\{ \begin{array}{c} xw+y^2 \\
x^8+yw+z^2 \end{array} \right\}$ & $19$  \\ 
$1;2,2$ & 1,2,4/8 & $X_{2,0}$ & $x^8+y^4+z^2$ & $21$  \\ 
$1;2,3$ & 1,2,3/7 & $S^\ast_{2,0}$ & $x^7+xz^2+y^2z$ & $20$  \\
$1;2,4$ & 1,2,3,4/5,6 & $L^\ast_{2,0}$ & $\left\{ \begin{array}{c} xw+yz \\
x^6+yw+z^2 \end{array} \right\}$ & $19$  \\ 
$1;3,3$ & 1,2,3,3/4,6 & $K'X_{2,0}$ & $\left\{ \begin{array}{c} xw+y^2 \\
x^6+z^2+w^2 \end{array} \right\}$ & $19$ \\ 
$1;2,2,2$ & 1,2,2/6 & $U^\ast_{2,0}$ & $x^6+y^3+z^3$ & $20$ \\ 
$1;2,2,3$ & 1,2,2,3/4,5 & $M^\ast_{2,0}$ & $\left\{ \begin{array}{c} xw+yz \\
(x^3+w)(y+z) \end{array} \right\}$ & $19$  \\ 
$1;2,2,2,2$ & 1,2,2,2/4,4 & $IT_{2,2,2,2}$ & $\left\{ \begin{array}{c} x^4+(y-z)w \\
ax^4+y(z-w) \end{array} \right\}$ & $19$  \\ 
$2;$ & 1,1,3/6 &  &  $x^6+y^6+z^2$ & $25$  \\
$2;2$ & 1,1,2/5 &  & $x^5+xz^2+y^5$ & $24$  \\ 
$2;3$ & 1,1,2,3/4,4 &  & $\left\{ \begin{array}{c} xw+y^4 \\
x^4+yw+z^2 \end{array} \right\}$ & $23$  \\ 
$2;2,2$ & 1,1,2,2/3,4 &  & $\left\{ \begin{array}{c} xw+yz \\
x^4+y^4+zw \end{array} \right\}$
& $23$ \\ 
$3;({\rm nh})$ & 1,1,1/4 &  & $x^4+y^4+z^4$ & $27$  \\
$3;({\rm h})$ & 1,1,1,2/2,4 &  & $\left\{ \begin{array}{c} xy+z^2 \\
x^4+y^4+z^4+w^2 \end{array} \right\}$ & $27$  \\ 
$3;2({\rm nh})$ & 1,1,1,2/3,3 &  & $\left\{ \begin{array}{c} xw+y^3 \\
x^3+yw+z^3 \end{array} \right\}$ & $26$ \\ 
$3;2({\rm h})$ & 1,1,1,2,2/2,3,3 &  & & $26$ \\ 
$4;({\rm nh})$ & 1,1,1,1/2,3 &  & $\left\{ \begin{array}{c} xw+yz \\
x^3+y^3+z^3+w^3 \end{array} \right\}$ & $29$  \\ 
$4;({\rm h})$ & $\begin{array}{c}1,1,1,1,2,2/ \\
2,2,2,3 \end{array}$ &  & & $29$ \\ 
$5;({\rm nh})$ & $\begin{array}{c}1,1,1,1,1/ \\
2,2,2 \end{array}$ &  & 3 quadrics in $\CC^5$ & $31$
\\
$5;({\rm h})$ & $\begin{array}{c}1,1,1,1,1,2,2,2/ \\
2,2,2,2,2,2 \end{array}$ &   & & $31$ \\ 
\hline
\end{tabular}}
\end{table}

This leads to the following classification. The cases $g=0$ and $r \leq 4$ are listed
in Table~\ref{table1}, the remaining cases in Table~\ref{table2}.   Here we
use the following notation.
We first list the orbit invariants. For $g
\geq 3$ we add (h) or (nh) to indicate whether $Z$ is hyperelliptic or
non-hyperelliptic respectively. In the second column we give the weights and the degrees
of the singularity. In the third column we indicate the name of the singularity
according to Arnold's
\cite{Arnold76} or Wall's notation \cite{Wall83, Wall84}, if it exists.  In
column 4 we list the equation(s) of the singularity. Here $a$ is a complex number with $a
\neq 0,1$. The cases $g=0$ and
$r\leq 4$  are Kodaira singularities in the sense of \cite{EW85}. They were already
considered in
\cite{Ebeling96}. The remaining cases with $g=0$ are still minimally elliptic  and
equations are given in
\cite{Ebeling86} and \cite{Wall91} respectively. Finally we indicate the Milnor number
$\mu$.

\addvspace{3mm}

\noindent {\bf Remark 2}\ A Fuchsian singularity with signature $\{5;\}$ is the cone over
a canonical curve of genus $5$ in $\PP^4$. As R.-O.~Buchweitz pointed out to me, this
is in general not an ICIS, but given by the Pfaffians of a $5 \times 5$ skew symmetric
matrix.

\addvspace{3mm}


Now let $(X,x)$ be a Fuchsian ICIS. All 
Fuchsian ICIS in $\CC^3$ and $\CC^4$ have equations such that one
of the Theorems~\ref{phi_M3}, \ref{phi_M4a}, and
\ref{phi_M4b} is applicable. The remaining singularities are 3 ICIS with
hyperelliptic $Z$ and the  singularity
$5;{\rm (nh)}$. There is the following relation between the singularities with the same
signature, but with hyperelliptic $Z$ on the one hand and non-hyperelliptic $Z$ on the
other hand. If the one with non-hyperelliptic $Z$ has 
weights $q_1, \ldots , q_n$ and degrees
$d_1, \ldots, d_{n-2}$, then the one with the same signature but hyperelliptic $Z$ has
weights $q_1, \ldots, q_n, q'_1, \ldots , q'_m$ and degrees $d_1, \ldots, d_{n-2}, q'_1,
\ldots , q'_m$ for some $m$. But this means that both singularities have 
the same Poincar\'{e} series and the same characteristic
polynomial of the monodromy and Theorem~\ref{phi_M4a} or
Theorem~\ref{phi_M4b} is also applicable in this case. Finally, for the singularity
$5;{\rm (nh)}$ we have
$\phi_M(t)=(1-t^2)^{16}/(1-t)$ and we define
$\phi_M^\flat(t):=\phi_M(t)/(1-t)^7=(1-t^2)^{16}/(1-t)^8$ and
$\tilde{\phi}_A(t):=(1-t^2)^{8}/(1-t)^{16}$. If $(X,x)$
is a hypersurface singularity, then we set
$\phi^\flat_M(t):=\phi_M(t)$. Therefore for all Fuchsian ICIS, 
polynomials
$\tilde{\phi}_A(t)$ and $\phi^\flat_M(t)$ are
defined and we obtain:

\begin{corollary}
If $(X,x)$ is a Fuchsian ICIS, then one has 
$$\tilde{\phi}_A^\ast(t) = \phi^\flat_M(t).$$
\end{corollary}

There is the following relation of this duality  with the mirror symmetry of
$K3$ surfaces. 

Let
$(X,x)$ be a normal surface singularity with good $\CC^\ast$-action. According to
\cite{Pinkham74}, the $\CC^\ast$-action on $X$ extends (at least formally) to a
$\CC^\ast$-action on a semi-universal deformation $p: \cal{X} \to S$ of $(X,x)$. Assume
that $(X,x)$ is smoothable. 
Then also the Milnor fibre ${\cal X}_t$  can be compactified in a
natural way to a surface  $\bar{\cal X}_t$ with the same cyclic quotient singularities
as $\bar{X}$ along a curve at infinity isomorphic to $\bar{X}_\infty$.  Denote by $Y$ the
minimal resolution of $\bar{\cal X}_t$. By 
\cite[6.13 Theorem]{Pinkham78}, $Y$ is a minimal $K3$ surface if and only if $(X,x)$ is
Fuchsian. 

Let $(X,x)$ be a Fuchsian ICIS. Let $Y^-$ be a tubular neighbourhood of 
the curve at infinity $Y_\infty \cong \tilde{X}_\infty$, and set
$Y^+:=Y\setminus{\rm int}(Y^-)$ and $\Sigma:=\partial Y^+ = \partial Y^- = Y^+ \cap
Y^-$. Then the Mayer-Vietoris sequence of the pair $(Y^+,Y^-)$ gives the following
exact sequence (we consider homology with integral coefficients):
$$0 \to H_2(\Sigma) \to H_2(Y^+) \oplus H_2(Y^-) \to H_2(Y) \to H_1(\Sigma) \to
H_1(Y^+) \oplus H_1(Y^-) \to 0$$
Now $M:=H_2(Y^+)$ is the Milnor lattice. The group $H_2(\Sigma)$ is the radical $M_0:=
{\rm rad}\, (M)$ of $M$ (cf.\ e.g.\ \cite{HM68}). The rank of $M_0$ is 
$\mu_0 =2g$. Let $H_1(\Sigma)_{{\rm tors}}$ denote the torsion subgroup of
$H_1(\Sigma)$. Then
$H_1(\Sigma)/H_1(\Sigma)_{{\rm tors}} \cong H_1(Y^-) \cong \ZZ^{\mu_0}$ 
and
$H_1(Y^+)=0$. This shows that we have a primitive embedding of the non-degenerate
lattice $M/M_0$ corresponding to the Milnor lattice into the homology lattice
$K:=H_2(Y) \cong (-E_8) \oplus (-E_8) \oplus U  \oplus U \oplus U$ of a $K3$ surface
(cf.\ \cite{Pinkham77b}).  By a result of V.~V.~Nikulin
\cite{Nikulin79}, this embedding is unique up to isometries. Moreover, the orthogonal
complement of $M/M_0$ in $K$ is the lattice $H_2(Y^-)=M_\infty$.

The duality between the
lattices $M/M_0$ and $M_\infty$ corresponds to the mirror symmetry of
$K3$-surfaces, see \cite{Dolgachev96}. In the case when $g=0$ and $r \leq 4$, the
lattice $M_\infty \oplus U$ is related to the Milnor lattice and 
$c_\infty$ to the monodromy
operator of another Fuchsian singularity $(X^\ast,x^\ast)$ with
$g=0$ and $r \leq 4$ and we
obtain Arnold's strange duality and its extension \cite{Ebeling96}. It is not known
to the author whether $M_\infty \oplus U$ and $c_\infty$ 
correspond to a singularity in the other cases. 

\section{Singular moonshine} \label{SM}
We shall now discuss relations to the Leech lattice.

For a polynomial
$$\phi(t)= \prod_{m|h} (1 - t^m)^{\chi_m}$$
where $\chi_m \in \ZZ$,
we use the symbolic
notation
$$\pi:= \prod_{m|h} m^{\chi_m}.$$
In the theory of finite groups, this symbol is known as a {\em Frame shape} (cf.\
\cite{CN79}). 

Let $(X,x)$ be a singularity which satisfies the assumptions of Theorem~\ref{phi_M3},
\ref{phi_M4a}, or \ref{phi_M4b}. Let $\pi_M$ and $\pi_M^\flat$ denote the Frame shapes
corresponding to the polynomials $\phi_M(t)$ and $\phi_M^\flat(t)$ respectively. If
$(X,x)$ is a  hypersurface singularity, then we put $\pi:=\pi_M$. Otherwise we define
$\pi:=\pi_M^\flat$ according to Theorem~\ref{phi_M4a} or \ref{phi_M4b}.

Now let $(X,x)$ be a Fuchsian ICIS. All 
Fuchsian ICIS in $\CC^3$ and $\CC^4$ have equations such that one
of the Theorems~\ref{phi_M3}, \ref{phi_M4a}, and
\ref{phi_M4b} is applicable. Above we also defined a polynomial $\phi^\flat_M$ in the
remaining cases and we set $\pi:=\pi_M^\flat$ in these cases as well. Then we make the
following observation.

\begin{theorem} \label{Leech}
Let $(X,x)$ be a Fuchsian ICIS. Then
the symbol $\pi \pi^\ast$ is a $24$-dimensional self-dual Frame
shape which is the Frame shape of an automorphism of the Leech lattice.
\end{theorem}

In \cite[Appendix~1]{Saito98} K.~Saito considers regular systems of weights of
rank 24 having negative exponents. They correspond to quasihomogeneous hypersurface
singularities in $\CC^3$ with $R \geq 1$ and $\mu=24$. He finds 11 cases with $R>1$.
It can be checked using the normal forms of quasihomogeneous functions in three
variables \cite[13.2]{Arnold74} that the list is complete up to one case which was
omitted. The complete list is in Table~\ref{table4.3a}. By a computer search we found 13
quasihomogeneous ICIS in $\CC^4$ with
$\mu=25$. They are listed in Table~\ref{table4.3b}. We use the same conventions as for
Table~\ref{table1} and \ref{table2}. All singularities have $\mu_0=g=0$,
$p_g=2$, and hence
$\mu_+=4$.
\begin{table}\centering
\caption{Quasihomogeneous ICIS in $\CC^3$ with $\mu=24$}
\label{table4.3a}
\begin{tabular}{|c|c|c|c|c|c|c|}
\hline
$g; \alpha_1, \ldots , \alpha_r$ & $b$ & $R$ & Weights & Name & Equation &  \\
\hline
$0;7,7,7$ & 1 & 4 & $3,7,7/21$ & $U_{24}$ & $x^7+y^3+z^3$ & (b,e)\\
$0;3,7,11$ & 1 & 4 &  $3,7,11/25$ & $S_{24}$ & $x^6y+xz^2+y^2z$ & (b) \\
$0;3,3,15$ & 1 & 4 & $3,11,15/33$ & $Q_{24}$ & $x^{11}+xz^2+y^3$ & (b,e)\\
$0;4,7,9$ & 1 & 5 & $4,7,9/25$ & $V^{\sharp}NC^1_{18}$ & $x^4z+xy^3+yz^2$ & (b)\\
$0;4,4,12$ & 1 & 5 & $4,7,12/28$ & $VNC^1_{18}$ & $x^7+xz^2+y^4$ & (b,e)\\
$0;2,9,9$ & 1 & 5 & $4,9,18/36$ & $W_{24}$ & $x^9+y^4+z^2$ & (b,e)\\
$0;2,4,14$ & 1 & 5 & $4,14,23/46$ & $Z_{24}$ & $x^8y+xy^3+z^2$ & (b,e)\\
$0;3,6,9$ & 1 & 7 & $6,9,11/33$ & $V'(Z_{12})^2$ & $x^4y+xy^3+z^3$ & (b,e)\\
$0;2,6,10$ & 1 & 7 & $6,10,23/46$ & $NC^1_{17}$ & $x^6y+xy^4+z^2$ & (b,e)\\
$0;2,3,13$ & 1 & 7 & $6,26,39/78$ & $E_{24}$ & $x^{13}+y^3+z^2$ & (b,e)\\
$0;2,5,7$ & 1 & 11 & $10,14,35/70$ & $NF^1_{(4)}$ & $x^7+y^5+z^2$ & (b,e)\\
$0;3,4,5$ & 1 & 13 & $12,15,20/60$ & $V'(W_{12})^2$ & $x^5+y^4+z^3$ & (b,e)\\
\hline
\end{tabular}
\caption{Quasihomogeneous ICIS in $\CC^4$ with $\mu=25$}
\label{table4.3b}
{\tabcolsep3pt
\begin{tabular}{|c|c|c|c|c|c|c|}
\hline
$g; \alpha_i$ & $b$ & $R$ & Weights & Name & Equations &  \\
\hline
$0;5,5,5,5$ & 2 & 3 & $2,5,5,5/10,10$ & $I^2_{1,0}$ & 
$\left\{ \begin{array}{c} x^5+(y-z)w\\
ax^5+y(z-w) \end{array} \right\}$ & (b,e) \\
$0;2,5,5,8$ & 2 & 3 & $2,5,5,8/10,13$ & $M_{2,0}$ & 
$\left\{ \begin{array}{c} xw+yz\\
(x^4+w)(y+z) \end{array} \right\}$ & (b)\\
$0;2,2,8,8$ & 2 & 3 & $2,5,8,8/10,16$ & $K'W_{2,0}$ & 
$\left\{ \begin{array}{c} xz+y^2\\
x^8+z^2+w^2 \end{array} \right\}$ & (b,e)\\
$0;2,2,5,11$ & 2 & 3 & $2,5,8,11/13,16$ & $L_{2,0}$ & 
$\left\{ \begin{array}{c} xw+yz\\
x^8+yw+z^2 \end{array} \right\}$ & (b)\\
$0;2,2,2,14$ & 2 & 3 & $2,8,11,14/16,22$ & $J'_{4,0}$ & 
$\left\{ \begin{array}{c} xw+y^2\\
x^{11}+yw+z^2 \end{array} \right\}$ & (b,e)\\
$0;3,3,5,7$ & 2 & 4 & $3,5,6,7/12,13$ &  & $\left\{ \begin{array}{c} yw+z^2\\
x^2w+xy^2+zw \end{array} \right\}$ & \\
$0;3,3,3,9$ & 2 & 4 & $3,5,6,9/12,15$ &  & $\left\{ \begin{array}{c} xw+z^2\\
x^5+y^3+zw \end{array} \right\}$ & (e)\\
$0;2,4,4,6$ & 2 & 5 & $4,6,7,8/14,16$ &  & $\left\{ \begin{array}{c} yw+z^2\\
x^4+xy^2+w^2 \end{array} \right\}$ & (e)\\
$0;2,3,4,7$ & 2 & 5 & $4,6,7,9/13,18$ &  & $\left\{ \begin{array}{c} xw+yz\\
xz^2+y^3+w^2 \end{array} \right\}$ & (b)\\
$0;2,2,4,8$ & 2 & 5 & $4,6,8,11/12,22$ &  &
$\left\{ \begin{array}{c} x^3+y^2+xz \\
xy^3+x^4y+yz^2+w^2 \end{array} \right\}$ & (b) \\
$0;2,3,3,4$ & 2 & 7 & $6,8,9,12/18,24$ &  & $\left\{ \begin{array}{c} xw+z^2\\
x^4+xz^2+y^3+w^2 \end{array} \right\}$ & \\
$0;2,2,3,5$ & 2 & 7 & $6,8,10,15/16,30$ &  & $\left\{ \begin{array}{c} xz+y^2\\
x^5+xy^3+z^3+w^2 \end{array} \right\}$ & (b,e)\\
$0;3,4,20$ & 1 & 11 & $9,12,16,20/32,36$ &  & $\left\{ \begin{array}{c} yw+z^2\\
x^4+y^3+zw \end{array} \right\}$ & \\
\hline 
\end{tabular}}
\end{table}

\begin{table}\centering
\caption{Self-dual Frame shapes of $\cdot 0$}
\label{table6}
\begin{tabular}{|c|c|l|}  \hline
ATL & Frame shape  & realized 
\\ \hline & & \\
1A & $2^{24}/1^{24}$  & (a) $A_1^{24}$ (c) $\tilde{D}_5^4$ (d) $\{ 5; \}$ 
\\ 
3A & $3^{12}/1^{12}$  & (a) $A_2^{12}$ (c) $\tilde{E}_6^3$ (d) $\{3;2\}$, $\{4;\}$ 
\\ 
3B & $2^66^6/1^63^6$  & (a) $D_4^6$ (d) $\delta1$ 
\\ 
4A & $4^8/1^8$  & (a) $A_3^8$ (d) $IT_{2,2,2,2}$, $\{2;3\}$, $\{2;22\}$, $\{3;\}$ 
\\ 
5A & $5^6/1^6$  & (a) $A_4^6$ (d) $M^\ast_{2,0}$, $\{2;2\}$ (e) $I^2_{1,0}$ 
\\ 
5B & $2^410^4 / 1^45^4$ & (a) $D_6^4$ (d) $J'_{2,0}$, $NA^1_{0,0}$ 
\\ 
6A & $3^46^4 / 1^42^4$  & (d) $I_{1,0}$ 
\\ 
6D & $2 \! \cdot \! 6^5/1^53$  & (a) $A_5^4D_4$ \\
   &                           & (d) $\alpha1^{(1)}$, $\alpha1^{(2)}$, 
$L^\ast_{2,0}$, $K'X_{2,0}$, $U^\ast_{2,0}$, $\{2;\}$
\\ 
7A & $7^4 / 1^4$   & (a) $A_6^4$ (d) $M_{1,0}$, $S^\ast_{2,0}$ (e) $U_{24}$
\\ 
7B & $2^314^3 / 1^37^3$   & (a) $D_8^3$ (d) $J'_{10}$, $Z_{1,0}$ (e) $J'_{4,0}$ 
\\ 
8C & $2^28^4 / 1^44^2$   & 
(a) $A_7^2D_5^2$ (d) $L_{1,0}$, $K'_{1,0}$, $VNA^1_{0,0}$, $J'_{3,0}$, $X_{2,0}$ \\
   &                     & (e) $K'W_{2,0}$  
\\ 
9A & $9^3 / 1^3$   & (a) $A_8^3$ (d) $M_{11}$, $U_{1,0}$, $Q_{3,0}$ (e) $\{0;3339\}$
\\ 
9C & $2^33^218^3 /1^36^29^3$  & (a) $D_{10}E_7^2$ (d) $J'_9$, $J_{3,0}$ 
\\ 
10A & $5^210^2 / 1^22^2$   & (b) $I^2_{1,0}$ 
\\ 
10E & $2 \! \cdot \! 10^3 / 1^35$  & (a) $A_9^2D_6$ (d) $L_{11}$, $K'_{11}$, $S_{1,0}$,
$Z_{2,0}$  
\\ 
11A & $2^222^2 / 1^211^2$   & (a) $D_{12}^2$ (b) $J'_{4,0}$, $\{0;2248\}$ (d)
$Z_{12}$  
\\ 
12A & $2^43^412^4 /1^44^46^4$   & (a) $E_6^4$ 
\\ 
12E & $4^212^2 / 1^23^2$   & (d) $U_{12}$ (e) $VNC^1_{18}$, $\{0;2446\}$ 
\\ 
12K & $2^23 \! \cdot \! 12^3 / 1^34 \! \cdot \! 6^2$  &
(a) $A_{11}D_7E_6$ \\
    &  & (d) $L_{10}$, $K'_{10}$, $J'_{11}$, $Q_{2,0}$, $W_{1,0}$, $J_{4,0}$
\\ 
13A & $13^2 / 1^2$   & (a) $A_{12}^2$ (b) $M_{2,0}$ (d) $S_{12}$ 
\\ 
15A & $2^33^35^330^3 / 1^36^310^315^3$   & (a) $E_8^3$ 
\\ 
15B & $3^215^2 / 1^25^2$  & (d) $Q_{12}$ (e) $Q_{24}$ 
\\ 
15D & $2 \! \cdot \! 6 \! \cdot \! 10 \! \cdot \! 30/1 
\! \cdot \! 3 \! \cdot \! 5 \! \cdot \! 15$ & (e) $NC^1_{17}$  
\\ 
15E & $2^23 \! \cdot \!5 \! \cdot \! 30^2 / 1^26 \! \cdot \! 10 \! \cdot \! 15^2$ 
& (a) $D_{16}E_8$ (b,e) $\{0;2235\}$, (d) $Z_{11}$, $E_{13}$ 
\\  
16B & $2 \! \cdot \! 16^2 / 1^28$  & (a) $A_{15}D_9$
(b) $K'W_{2,0}$, $L_{2,0}$ (d) $S_{11}$, $W_{13}$
\\ 
18A & $9 \! \cdot \! 18 / 1 \! \cdot \! 2$ & (e) $W_{24}$, $V'(Z_{12})^2$ 
\\ 
18B & $2 \! \cdot \! 3 \! \cdot \! 18^2 / 1^26 \! \cdot \! 9$   & 
(a) $A_{17}E_7$ (b) $\{0;2347\}$ (d) $Q_{11}$, $Z_{13}$  
\\ 
20A & $2^25^220^2 / 1^24^210^2$   & (d) $W_{12}$ 
\\ 
21A & $2^23^27^242^2 / 1^26^214^221^2$  & (d) $E_{12}$ 
\\ 
21B & $7 \! \cdot \! 21 / 1 \! \cdot \! 3$ & (b) $U_{24}$ 
\\ 
23A & $2 \! \cdot \! 46/1 \! \cdot \! 23$ & (a) $D_{24}$ (b) $Z_{24}$, $NC^1_{17}$  
\\ 
24B & $2 \! \cdot \! 3^24 \! \cdot \! 24^2 / 
1^26 \! \cdot \! 8^2 12$  & (d) $Q_{10}$, $E_{14}$ 
\\ 
25A & $25/1$ & (a) $A_{24}$ (b) $S_{24}$, $V^{\sharp}NC^1_{18}$ 
\\ 
28A & $4 \! \cdot \! 28 / 1 \! \cdot \! 7$ & (b) $VNC^1_{18}$ (e) $Z_{24}$ 
\\ 
33A & $3 \! \cdot \! 33/1 \! \cdot \! 11$ & (b) $Q_{24}$, $V'(Z_{12})^2$ 
\\ 
35A & $ 2 \! \cdot \! 5 \! \cdot \! 7 \! \cdot \! 70 / 
1 \! \cdot \! 10 \! \cdot \! 14 \! \cdot \! 35$ & (b,e) $NF^1_{(4)}$  
\\ 
36A & $2 \! \cdot \! 9 \! \cdot \! 36/ 
1 \! \cdot \! 4 \! \cdot \! 18$ & (b) $W_{24}$ 
\\ 
39A & $2 \! \cdot \! 3 \! \cdot \! 13 \! \cdot \! 78 / 
1 \! \cdot \! 6 \! \cdot \! 26 \! \cdot \! 39$ & (b,e) $E_{24}$ 
\\ 
60A & $3 \! \cdot \! 4 \! \cdot \! 5 \! \cdot \! 60 / 
1 \! \cdot \! 12 \! \cdot \! 15 \! \cdot \! 20$ & (b,e) $V'(W_{12})^2$ 
\\ 
& & \\ \hline 
\end{tabular}
\end{table}

Saito already observed that the Frame shapes $\pi$ of the hypersurface
singularities are self-dual and appear as Frame shapes of automorphisms of the Leech
lattice. For 9 of the ICIS in $\CC^4$ with $\mu=25$ 
	we can apply Theorem~\ref{phi_M4a} or \ref{phi_M4b} and define
$\pi:=\pi_M^\flat$. This symbol
is also a self-dual Frame shape of dimension 24 which appears as the Frame shape of an
automorphism of the Leech lattice. In the remaining cases either Theorems~\ref{phi_M4a}
and \ref{phi_M4b} are not applicable and so
$\pi_M^\flat$ is not defined or $\pi_M^\flat$ is not self-dual.

In Table~\ref{table6} we have listed the 39 self-dual Frame shapes of the
automorphism group $\cdot 0$ of the Leech lattice. We use the ATLAS notation
\cite{ATLAS} for the conjugacy classes. We consider the following five
constructions. For a Frame shape
$\pi= \prod_{m|h}m^{\chi_m}$ we call the minimal 
$h$ the order of $\pi$ and $\deg \, \pi := \sum_{m|h}m\chi_m$ the degree of $\pi$.

(a) Consider any combination (direct sum)
$\pi= \pi_1 \cdots \pi_s $
of the (self-dual) Frame shapes of the Coxeter elements of the root systems of type
$A_l$, $D_l$, $E_6$, $E_7$, or $E_8$ such that the orders of the
$\pi_i$ are the same and the degree of $\pi$ is equal to $24$. (There are 23
such combinations; they correspond to the 23 Niemeier lattices, see e.g.\
\cite[Proposition~3.4]{Ebeling94}.)

(b) Consider the (self-dual) symbol $\pi$ of a singularity
of Table~\ref{table4.3a} or \ref{table4.3b} respectively for which (b) is indicated in
the last column.

(c) Consider the symbol $\pi \pi^\ast$ where 
$\pi= \pi_1 \cdots \pi_s $
is any combination (direct sum) of Frame shapes of the simply elliptic singularities (cf.\
\cite{Ebeling01}) such that the orders of the
$\pi_i$ are the same and the degree of $\pi$ is equal to $24$. (There are
only two such combinations, namely $\tilde{D}_5^4$ and $\tilde{E}_6^3$.)

(d) Consider the symbol $\pi \pi^\ast$ where $\pi$ is the symbol of a Fuchsian ICIS
according to Theorem~\ref{Leech}.

(e) Consider the symbol $\pi \pi^\ast$ where $\pi$ is the symbol corresponding to the
polynomial $\psi_A(t)$ of a singularity of Table~\ref{table4.3a} or
\ref{table4.3b} for which (e) is indicated in the last column (cf.\ Remark~1).

By these constructions we get all of the self-dual Frame shapes of $\cdot 0$ (cf.\ also
\cite[Appendix~1]{Saito98}, where 4 cases do not appear). The
different realizations are indicated in Table~\ref{table6}. 

To a Frame shape $\pi = \prod_{m | h} m^{\chi_m}$
one can associate a modular function \cite{Kondo85}. Let
$$\eta(\tau) = q^{1/24} \prod^\infty_{n=1} (1-q^n), \quad q=e^{2\pi i \tau},
\tau \in \HH,$$
be the Dedekind eta function. Then define
$$\eta_\pi(\tau) = \prod_{m | h} \eta(m\tau)^{\chi_m}.$$
Let $\pi$ be a self-dual Frame shape of $\cdot 0$. By \cite{Kondo85},
$\eta_\pi$ is a modular function for a discrete subgroup $\Gamma'$ of 
$SL(2,\RR)$ containing $\Gamma_0(h)$. The genus of $\Gamma'$
is zero and $\eta_\pi$ is a generator of the function field of $\Gamma'$. The
groups $\Gamma'$ corresponding to the Frame shapes of Table~\ref{table6} are listed in
\cite{CN79}.

\end{document}